\documentclass[a4paper, 11pt, final]{article}

\usepackage{amsfonts,enumerate}
\usepackage{enumerate}

\newtheorem{proposition}{Proposition}[section]
\newtheorem{theorem}[proposition]{Theorem}
\newtheorem{lemma}[proposition]{Lemma}
\newtheorem{corollary}[proposition]{Corollary}

\newenvironment{proof}{\smallskip\noindent\emph{Proof.}\hspace{1pt}}%
{\hspace{-5pt}{\nobreak\quad\nobreak\hfill\nobreak$\square$\vspace{8pt}%
\par}\smallskip\goodbreak}

\newenvironment{proofof}[1]{\smallskip\noindent\emph{Proof of #1.}%
\hspace{1pt}}{\hspace{-5pt}{\nobreak\quad\nobreak\hfill\nobreak%
$\square$\vspace{8pt}\par}\smallskip\goodbreak}

\newcommand{\Section}[1]{\section{#1}\setcounter{equation}{0}}

\newcommand{\Id}{\mathinner{\mathrm{Id}}}

\newcommand{\C}[1]{\mathbf{C^{#1}}}
\newcommand{\PC}{\mathbf{PC}}
\newcommand{\Cc}[1]{\mathbf{C_c^{#1}}}

\newcommand{\modulo}[1]{{\left|#1\right|}}
\newcommand{\norma}[1]{{\left\|#1\right\|}}
\newcommand{\Ref}[1]{{\rm(\ref{#1})}}
\newcommand{\reali}{{\mathbb{R}}}
\newcommand{\naturali}{{\mathbb{N}}}
\newcommand{\tv}{\mathrm{TV}}
\newcommand{\BV}{\mathbf{BV}}
\newcommand{\Lip}{\mathinner\mathbf{Lip}}

\renewcommand{\O}{\mathinner{\mathcal{O}(1)}}
\renewcommand{\epsilon}{\varepsilon}
\renewcommand{\phi}{\varphi}
\renewcommand{\theta}{\vartheta}
\renewcommand{\L}[1]{\mathbf{L^#1}}

\title{Hyperbolic Balance Laws with a Non Local Source}

\author{Rinaldo M.~Colombo \\ \small Dipartimento di Matematica \\
  \small Universit\`a degli Studi di Brescia \\ \small Via Branze, 38
  \\ \small 25123 Brescia, Italy \\ \texttt{Rinaldo.Colombo@unibs.it}\\
  \and Graziano Guerra \\ \small Dip.~di Matematica e Applicazioni \\
  \small Universit\`a di Milano -- Bicocca \\ \small Via Bicocca degli
  Arcimboldi, 8 \\ \small 20126 Milano, Italy \\
  \texttt{Graziano.Guerra@unimib.it}}

\begin{document}

\maketitle

\begin{abstract}

  This paper is devoted to hyperbolic systems of balance laws with non
  local source terms. The existence, uniqueness and Lipschitz
  dependence proved here comprise previous results in the literature
  and can be applied to physical models, such as Euler system for a
  radiating gas and Rosenau regularization of the Chapman-Enskog
  expansion.

  \medskip

  \noindent\textit{2000~Mathematics Subject Classification:} 35L65,
  76N15

  \medskip

  \noindent\textit{Key words and phrases:} Systems of Balance Laws,
  Radiating Gas.
\end{abstract}

\Section{Introduction}
\label{sec:Intro}

This paper is devoted to systems of conservation laws with non local
sources, i.e.~to equations of the form
\begin{equation}
  \label{eq:Main}
  \partial_t u + \partial_x f(u) = G(u)
\end{equation}
where $f$ is the flow of a nonlinear hyperbolic system of conservation
laws and $G \colon \L1 \mapsto \L1$ is a (possibly) \emph{non local}
operator.  As examples, we consider below the case $G(u) = g(u) + Q *
u$ that enters a classical radiating gas model,
see~\cite{VincentiKruger}, as well as Rosenau regularization of
Chapman-Enskog expansion of the Boltzmann equation, see~\cite{Rosenau,
  SchochetTadmor}. We establish that~\Ref{eq:Main} is well posed in
$\L1$, locally in time, for data having sufficiently small total
variation. To this aim, we require on~\Ref{eq:Main} those assumptions
that separately guarantee the well posedness of the convective part
\begin{equation}
  \label{eq:HCL}
  \partial_t u + \partial_x f(u) = 0
\end{equation}
and of the source part
\begin{equation}
  \label{eq:Source}
  \partial_t u = G(u) \,.
\end{equation}
These two equations generate two semigroups of solutions, say $S$ and
$\Sigma$. To obtain our results we exploit the techniques
in~\cite{AmadoriGuerra2002, ColomboCorli5}, essentially based on the
\emph{fractional step} algorithm, see~\cite{ColomboCorli5,
  CrandallMajda, DafermosHsiao, Trotter}. Its core idea is to get a
solution of the original equation as a limit of approximations
obtained suitably merging $S$ and $\Sigma$.

On the two semigroups we require the following two key conditions
(see~\cite[\textbf{(S2)} and~\textbf{(C*)}]{ColomboCorli5}):
\begin{enumerate}[{\it i)}]
\item A Gr\"onwall type estimate with respect to a suitable metric
  $d(\cdot,\cdot)$: i.e.~for a positive $C$
  \begin{equation}
    \label{Gronwall}
    \begin{array}{l}
      d(S_tu,S_tv) \leq e^{Ct}d(u,v)
      \\
      d(\Sigma_t u,\Sigma_t v) \leq e^{Ct} d(u,v)
    \end{array}
    \mbox{ for all }t \geq 0.
  \end{equation}
\item A commutativity relation
  \begin{equation}
    \label{commutation}
    d(\Sigma_t S_t u, S_t \Sigma_t u) \leq K \, t^2
    \mbox{ as } t \to 0.
  \end{equation}
\end{enumerate}
\noindent The former assumption is used to prove the uniformly
continuous dependence of the approximations from the initial data,
also called \emph{stability condition} in the framework of Lie-Trotter
formula, see~\cite[Corollary~5.8, Chapter~3]{EngelNagel}. The latter
condition yields the convergence of the approximations and ensures the
uniqueness of the limit.

Assume that{\it~i)} holds and $d$ is a reasonable (in the sense of
Proposition~\ref{prop:Guerra}) metric equivalent to the $\L1$
distance. Then, the invariance under the hyperbolic rescaling $(t,x)
\to (\lambda t, \lambda x)$ of solutions to system of conservation
laws, implies that $C=0$, see~Proposition~\ref{prop:Guerra}. Hence, to
apply the operator splitting techniques, we need a contractive metric
for the conservation law~\Ref{eq:HCL}. This role is naturally played
by the well known functional $\Phi$ in~\cite{LiuYang1, LiuYang3}.
Note, however, that this functional is \emph{not} a metric, for it may
lack to satisfy the triangle inequality.  The proof, then, consists in
showing that the semigroup generated by the source
part~\Ref{eq:Source} satisfies~\Ref{Gronwall} with respect to Liu \&
Yang functional and commutes with the semigroup generated by the
conservation law in the sense of~\Ref{commutation}.

\medskip

More precisely, let $\Omega$ be a open subset of $\reali^n$ with $0
\in \Omega$. For all positive $\delta$, define
\begin{displaymath}
  \mathcal{U}_\delta
  =
  \left\{
    u \in \L1 (\reali;\Omega) \colon \tv(u) \leq \delta
  \right\} \,.
\end{displaymath}
As a general reference on conservation laws we refer
to~\cite{BressanLectureNotes}. On the convective and on the source
parts we assume throughout that
\begin{description}
\item[(F)] $f \in \C4 ( \Omega;\reali^n)$ is strictly hyperbolic and
  each characteristic field is either genuinely nonlinear or linearly
  degenerate.
\item[(G)] For a positive $\delta_o$, $G \colon \mathcal{U}_{\delta_o}
  \mapsto \L1(\reali,\reali^n)$ is such that for suitable positive
  $L_1,L_2,L_3$
  \begin{displaymath}
    \begin{array}{l@{\qquad}rcl}
      \forall\, u,w \in \mathcal{U}_{\delta_o} &
      \norma{G(u) - G(w)}_{\L1} & \leq & L_1 \cdot \norma{u-w}_{\L1}
      \\
      \forall\, u \in \mathcal{U}_{\delta_o} &
      \tv \left(G(u) \right) & \leq & L_2 \cdot \tv(u) + L_3\,.
    \end{array}
  \end{displaymath}
\end{description}

\noindent Note that~\textbf{(F)}, respectively~\textbf{(G)}, ensures
the local in time well posedness of~\Ref{eq:HCL},
respectively~\Ref{eq:Source}. A class of functions
satisfying~\textbf{(G)} is provided by the following proposition.

\begin{proposition}
  \label{prop:G}
  Let $g,h \colon \Omega \mapsto \reali^n$ be locally Lipschitz and $Q
  \in \L1 (\reali;\reali^{n\times n})$. Then, the operator $G(u) =
  g(u) + Q*h(u)$ satisfies~\textbf{(G)} with $L_3=0$.
\end{proposition}

\noindent The proof is deferred to Section~\ref{sec:Tech}. We are now
ready to state the main result of this work.

\begin{theorem}
  \label{thm:main}
  Let $f$ satisfy~\textbf{(F)} and $G$ satisfy~\textbf{(G)}. Then,
  there exist positive $T$, $\tilde\delta$, $\mathcal{L}$, closed
  domains $\mathcal{D}_t$ and processes
  \begin{displaymath}
    F_t \colon \mathcal{D}_{T-t} \mapsto \mathcal{D}_T
    \qquad \forall\,t \in [0,T]
  \end{displaymath}
  with the properties:
  \begin{enumerate}[(1)]
  \item for $t, s \in [0,T]$ with $t < s$, $\mathcal{U}_{\tilde
      \delta} \subseteq \mathcal{D}_{t} \subseteq \mathcal{D}_{s}
    \subseteq \mathcal{U}_{\delta_o}$;
  \item \label{it:thm:semi} for $u$ in $\mathcal{D}_T$, $F_0 u = u$;
    for $t,s \in [0,T]$ with $t + s \in [0,T]$, $F_{s} \mathcal{D}_{t}
    \subseteq \mathcal{D}_{t + s}$ and for $u \in \mathcal{D}_{T - t -
      s}$, $F_{t} F_{s} u = F_{t + s} u$;
  \item \label{it:3} for $\bar t \in [0,T]$ and $u \in
    \mathcal{D}_{\bar t}$, the map $t \mapsto F_{t} u$ is a weak
    entropy solution to~\Ref{eq:Main} for $t \in [0, T-\bar t]$;
  \item \label{it:tg}if $S$ is the SRS generated by~\Ref{eq:HCL}, then
    for $\bar t \in \left[0,T\right[$ and $u \in \mathcal{D}_{\bar
      t}$,
    \begin{displaymath}
      \lim_{t \to 0}
      \frac{1}{t} \,
      \norma{F_t u - \left( S_t u + t \, G(u) \right) }_{\L1} =0 \,;
    \end{displaymath}
  \item for $t,s \in [0,T]$, $u,w \in \mathcal{D}_{T-t}$ and $s<t$,
    then
    \begin{equation}
      \label{eq:error}
      \begin{array}{rcl}
        \displaystyle
        \norma{F_t u - F_t w}_{\L1} 
        & \leq &
        \displaystyle 
        \mathcal{L}_{\phantom{\vert}} \cdot \norma{u-w}_{\L1}
        \\
        \displaystyle
        \norma{F_t u - F_s u}_{\L1} 
        & \leq &
        \displaystyle
        \mathcal{L} \cdot \left( 1 + \norma{u}_{\L1} \right)^{\phantom{\vert}}
        \cdot \modulo{t-s}\,;
      \end{array}
    \end{equation}
  \item \label{it:thm:main} for $\bar t \in \left[0,T\right[$, $u \in
    \mathcal{D}_{\bar t}$ and $\tau \in [0, T-\bar t\,]$ the map
    $u(\tau) = F_\tau u$ satisfies
    \begin{enumerate}[(\ref{it:thm:main}a)]
    \item \label{it:thm:main1} for $\xi \in \reali$, $\displaystyle
      \lim_{\theta\to 0+} \frac{1}{\theta} \, \int_{\xi - \theta
        \hat\lambda}^{\xi+ \theta \hat\lambda} \norma{ \left(F_\theta
          u(\tau)\right)(x) - U^\sharp_{(u(\tau),\xi)}(\theta,x) } \,
      dx =0$,
    \item \label{it:thm:main2} there exists a positive $C$ such that
      for all $a,b$ with $-\infty \leq a < \xi < b \leq +\infty$
      \begin{eqnarray*}
        & &
        \limsup_{\theta \to 0+}
        \frac{1}{\theta} \,
        \int_{a + \theta\hat\lambda}^{b - \theta \hat\lambda} 
        \norma{\left(F_\theta u(\tau)\right)(x) - 
          U^\flat_{(u(\tau),\xi)} (\theta,x)} \, dx
        \quad \leq
        \\
        & \leq &
        C \cdot \left( \tv \left( u(\tau);\left]a,b\right[\right) \right)^2 \,,
      \end{eqnarray*}
    \end{enumerate}
    \noindent where $U^\sharp_{(u(\tau),\xi)}$ solves~\Ref{eq:Sharp}
    and $U^\flat_{(u(\tau),\xi)}$ solves~\Ref{eq:Flat};
  \item \label{it:thm2:main} for $\bar t \in \left[0, T \right[$, if a
    Lipschitz map $u \colon [0,T-\bar t] \mapsto \mathcal{D}_T$ is
    such that $u(t) \in \mathcal{D}_{\bar t+t}$ and for $\tau \in [0,
    T-\bar t]$
    \begin{enumerate}[(\ref{it:thm2:main}a)]
    \item \label{it:thm2:main1} for $\xi \in \reali$, $\displaystyle
      \lim_{\theta\to 0+} \frac{1}{\theta} \, \int_{\xi - \theta
        \hat\lambda}^{\xi+ \theta \hat\lambda} \norma{ \left(F_\theta
          u(\tau)\right)(x) - U^\sharp_{(u(\tau),\xi)}(\theta,x) } \,
      dx =0$,
    \item \label{it:thm2:main2} there exists a finite measure
      $\mu_\tau$ such that for all $a,b$ with $-\infty \leq a < \xi <
      b \leq +\infty$
      \begin{displaymath}
        \!\!\!\!
        \limsup_{\theta \to 0+}
        \frac{1}{\theta} \,
        \int_{a + \theta\hat\lambda}^{b - \theta \hat\lambda} 
        \norma{\left(F_\theta u(\tau)\right)(x) - 
          U^\flat_{(u(\tau),\xi)} (\theta,x)} \, dx
        \leq \left( \mu_\tau \left( \left]a,b \right[ \right) \right)^2
      \end{displaymath}
    \end{enumerate}
    \noindent where $U^\sharp_{(u(\tau),\xi)}$ solves~\Ref{eq:Sharp}
    and $U^\flat_{(u(\tau),\xi)}$ solves~\Ref{eq:Flat}, then $u(t) =
    F_t u(0)$.
  \end{enumerate}
  \noindent Moreover, if $f_1,f_2$ both satisfy~\textbf{(F)} and
  $G_1,G_2$ both satisfy~\textbf{(G)}, then, denoting by $F^i$ the
  process generated by $f_i$ and $G_i$, for all $t \in [0,T]$ and $u
  \in \mathcal{D}_0$
  \begin{equation}
    \label{eq:Uffa}
    \begin{array}{rcl}
      \displaystyle
      \norma{F^1_t u - F^2_t u}_{\L1}
      & \leq &
      \displaystyle
      \mathcal{L} \cdot \norma{Df_1-Df_2}_{\C0(\Omega,\reali^{n\times n})}
      \cdot t
      \\
      & &
      \displaystyle
      + \, 
      \mathcal{L} \cdot
      \norma{G_1-G_2}_{\C0(\mathcal{U}_{\delta_o};\L1(\reali;\reali^n))}
      \cdot t \,.
    \end{array}
  \end{equation}
\end{theorem}

For the definition and properties of the SRS, refer
to~\cite{BressanLectureNotes}. Point~\ref{it:tg}.~characterizes the
tangent vector to $t\mapsto F_t u$ in the sense
of~\cite[\S~5]{BressanLarnas}. It is through this characterization
that the integral inequalities~\ref{it:thm:main}
and~\ref{it:thm2:main} are proved.  The proof of
Theorem~\ref{thm:main} follows from the results presented in
Section~\ref{sec:Tech} below.

The first part of next section is devoted to the application of the
above result to the Euler system for a radiating gas and to Rosenau
regularization of the Chapman-Enskog expansion. The framework of local
sources is recovered in the subsequent paragraph and, finally, we
quickly comprise also the case of a non autonomous source.

Remark that in the estimate~\Ref{eq:error} the presence of the term
$\norma{u}_{\L1}$ is mandatory, as the example $\partial_t u = u$
shows. Indeed, the domains $\mathcal{U}_\delta$ is unbounded in $\L1$.
Moreover, note that the analogous estimate in~\cite{AmadoriGuerra2002,
  ColomboCorliJHDE, CrastaPiccoli} should be understood with a time
Lipschitz constant dependent on the $\L1$ norm of the initial datum.

We stress that the estimate~(\ref{it:thm:main2}) is sharper
than~\cite[formula~(5.18)]{AmadoriGosseGuerra} thanks to the finite
total variation of the source term, ensured by~\textbf{(G)}.

\Section{Applications and Extensions}
\label{sec:Appl}

\subsection{Euler System for a Radiating Gas}

The following model for a radiating polytropic gas was considered
in~\cite[Chapter~XXII, \S~6]{VincentiKruger}, see
also~\cite[formula~(1.2)]{LattanzioMarcati}:
\begin{displaymath}
  \left\{
    \begin{array}{l}
      \partial_t \rho + 
      \partial_x \left( \rho \, v \right)
      = 0
      \\
      \partial_t \left( \rho\, v\right) + 
      \partial_x \left( \rho\, v^2 + p \right)
      = 0
      \\
      \partial_t \left( \rho \, e + \frac{1}{2}\rho \, v^2 \right) +
      \partial_x \left( 
        v \left( \rho \, e + \frac{1}{2}\rho \, v^2 +p\right) +q
      \right)
      = 0
      \\
      -\partial_{xx}^2 q + a\, q+ b \, \partial_x \theta^4
      = 0
    \end{array}
  \right.
\end{displaymath}
Here, as usual, $\rho$ is the gas density, $v$ its speed, $e$ the
internal energy, $p$ the pressure, $\theta = e / c_v$ the temperature
and $q$ is the radiative heat flux. The system is closed by means of
the equation of state and specifying the values of the characteristic
constants $a$ and $b$.

Solving the latter equation in $q$ we have $q = -\frac{b}{\sqrt{a}} \,
Q_a * \left(\frac{d~}{dx}\theta^4\right)$, where $Q_a (x) =
\frac{1}{2}\, \exp \left( -\sqrt{a}\, \modulo{x}\right)$ and we are
lead to consider the system
\begin{equation}
  \label{eq:Radiating}
  \left\{
    \begin{array}{l}
      \partial_t \rho + 
      \partial_x \left( \rho \, v \right)_{\vphantom{\big\vert}}
      = 0
      \\
      \partial_t \left( \rho\, v\right) + 
      \partial_x \left( \rho\, v^2 + p \right)^{\vphantom{\big\vert}}
      = 0
      \\
      \partial_t \left( \rho \, e + \frac{1}{2}\rho \, v^2 \right) +
      \partial_x \left( 
        v \left( \rho \, e + \frac{1}{2}\rho \, v^2 +p\right)
      \right)
      = b \, \left( - \theta^4 + \sqrt{a}\, Q_a * \theta^4 \right). \!\!\! \!\!\! \!\!\! \!\!\! \!\!\! \!\!\! \!\!\! \!\!\!
    \end{array}
  \right.
\end{equation}
It is well known that Euler system satisfies~\textbf{(F)}.
Condition~\textbf{(G)} holds by Proposition~\ref{prop:G}. Hence,
Theorem~\ref{thm:main} applies and we obtain the local in time well
posedness of~\Ref{eq:Radiating}. Note that this result also ensures
the local Lipschitz dependence of the solutions to~\Ref{eq:Radiating}
from the parameters $a$ and $b$.

\subsection{Rosenau Regularization of the Chapman-Enskog Expansion}

In his classical work~\cite{Rosenau}, Rosenau proposed a system of
balance laws that provides a regularized version of the Chapman-Enskog
expansion for hydrodynamics in a linearized framework. The 1D version
is the following:
\begin{displaymath}
  \left\{
    \begin{array}{l}
      \partial_t \rho + \partial_x v =0
      \\
      \partial_t v + \partial_x p = \mu_* * \partial^2_{xx} v
      \\
      \partial_t \left( \frac{3}{2}\theta\right) + \partial_x v = 
      \lambda_* * \partial^2_{xx}\theta
    \end{array}
  \right.
\end{displaymath}
where $\rho$ is the fluid density, $v$ is its speed and $\theta$ is
the temperature. $\mu_*$, respectively $\lambda_*$, is a convolution
kernel related to viscosity, respectively to thermal conductivity.
This linear system motivated analytical results, see for
instance~\cite{JinSlemrod, LiuTadmor, SchochetTadmor}, mostly related
to the quasilinear scalar equation
\begin{displaymath}
  \partial_t u + \partial_x \left( \frac{1}{2}u^2\right) = -u + Q*u
\end{displaymath}
since the source term $-u + Q*u$ is equal to $Q*\partial^2_{xx}u$,
provided $Q(x) = \frac{1}{2} \exp \left( -\modulo{x} \right)$.
Therefore, it is natural to consider the following Euler system with
Rosenau-type sources
\begin{equation}
  \label{eq:Rosenau}
  \left\{
    \begin{array}{l}
      \partial_t \rho + 
      \partial_x \left( \rho \, v \right)
      = 0
      \\
      \partial_t \left( \rho\, v\right) + 
      \partial_x \left( \rho\, v^2 + p \right)
      = \mu_* * \partial^2_{xx} v
      \\
      \partial_t \left( \rho \, e + \frac{1}{2}\rho \, v^2 \right) +
      \partial_x \left( 
        v \left( \rho \, e + \frac{1}{2}\rho \, v^2 +p\right) +q
      \right)
      = \lambda_* * \partial^2_{xx}\theta \,.
    \end{array}
  \right.  
\end{equation}
Rosenau kernels, see~\cite[formul\ae~(4a) and~(6)]{Rosenau} read
\begin{displaymath}
  \mu_* (x)
  =
  \frac{\mu}{2\, m\, \epsilon} \, \exp \left( -\modulo{x}/\epsilon\right)
  \quad \mbox{ and }\quad
  \lambda_* (x)
  =
  \frac{\lambda}{2\,s\,\epsilon} \, \exp \left( -\modulo{x}/\epsilon\right)
\end{displaymath}
for suitable positive parameters $\mu,\lambda,m,s,\epsilon$. With the
above choices, the sources in the last two equations
in~\Ref{eq:Rosenau} can be rewritten as
\begin{displaymath}
  \mu_* * \partial^2_{xx} v 
  =
  \frac{1}{\epsilon^2} \left( -\frac{\mu}{m} \, v + \mu_* * v \right)
  \quad \mbox{ and } \quad
  \lambda_* * \partial^2_{xx} \theta
  =
  \frac{1}{\epsilon^2} 
  \left( -\frac{\lambda}{s} \, \theta + \lambda_* * \theta \right) \,.
\end{displaymath}
By~Proposition~\ref{prop:G}, system~\Ref{eq:Rosenau} falls within the
scope of Theorem~\ref{thm:main}. Thus, we prove the local in time well
posedness of~\Ref{eq:Rosenau} as well as the local Lipschitz
dependence of the solutions to~\Ref{eq:Radiating} from the parameters
$\mu,\lambda,m,s,\epsilon$.

\subsection{Local Inhomogeneous Source}

Theorem~\ref{thm:main} can be applied also in the standard case of a
\emph{local} source. Indeed, it is immediate to see that~\textbf{(G)}
is implied by the following conditions~\textbf{(g1)}
and~\textbf{(g2)}.

\begin{proposition}
  Let $g \colon \reali \times \Omega \mapsto \reali^n$ be such that
  \begin{enumerate}[{\bf(g1)}]
  \item there exists an $L_1 > 0$ such that for $u_1,u_2 \in \Omega$,
    $\norma{g(x,u_2) - g(x,u_1)} \leq L_1 \cdot \norma{u_2 - u_1}$;
  \item there exists a finite measure $\mu$ on $\reali$ such that for
    $u \in \Omega$ and $x_1,x_2 \in \reali$ with $x_1 < x_2$,
    $\norma{g(x_2,u) - g(x_1,u)} \leq \mu\left( \left] x_1, x_2
      \right] \right)$.
  \end{enumerate}
  \noindent and assume that $f$ satisfies~\textbf{(F)}. Then, setting
  $\left( G(u)\right)(x) = g(x,u)$, Theorem~\ref{thm:main} applies.
\end{proposition}

Note that the integral estimates~\ref{it:thm:main}
and~\ref{it:thm2:main} in Theorem~\ref{thm:main} ensure that the
solution constructed here coincide with those
in~\cite{AmadoriGosseGuerra}. Similarly, the
characterization~\ref{it:tg} of the tangent vector imply that the
present solutions coincide with those in~\cite{AmadoriGuerra2002}.

\subsection{The Non Autonomous Case}

Theorem~\ref{thm:main} can be extended to the non autonomous balance
law
\begin{displaymath}
  \partial_t u + \partial_x f(u) = G(t,u)
\end{displaymath}
provided $f$ satisfies~\textbf{(F)}, $G \colon [0, T_o] \times \L1
\mapsto \L1$ satisfies
\begin{description}
\item[(G')] For positive $\delta_o,T_o$, the map $G \colon [0,T_o]
  \times \mathcal{U}_{\delta_o} \mapsto \L1(\reali,\reali^n)$ admits
  suitable positive $L_1,L_2,L_3$ such that for all $u,w \in
  \mathcal{U}_{\delta_o}$ and for all $t,s \in [0,T_o]$
  \begin{displaymath}
    \begin{array}{rcl}
      \displaystyle
      \norma{G(t,u) - G(s,w)}_{\L1} 
      & \leq & 
      \displaystyle
      L_1 \cdot \left(\norma{u-w}_{\L1} + \modulo{t-s} \right)
      \\
      \displaystyle
      \tv \left(G(t,u) \right) 
      & \leq & 
      \displaystyle
      L_2 \cdot \tv(u) + L_3\,.
    \end{array}
  \end{displaymath}
\end{description}

\noindent Indeed, let $\hat \lambda$ be an upper bound for all moduli
of characteristic speeds, i.e.~$\hat \lambda > \sup_{\norma{u}\leq
  \delta_o} \max_{i=1,\ldots,n} \modulo{\lambda_i(u)}$, and define
\begin{displaymath}
  \tilde f(u,w) = \left( f(u),\hat\lambda\, w\right)
  \qquad
  \tilde G (u,w) = 
  \left( G\left( \textstyle\int_{\reali }w,u \right), \chi_{\strut[0,1]} \right)
\end{displaymath}
here, $\chi_{\strut[0,1]}$ is the characteristic function of the real
interval $[0,1]$. Then, $\tilde f$ satisfies~\textbf{(F)} and $\tilde
G$ satisfies~\textbf{(G)}, so that Theorem~\ref{thm:main} applies and
the balance law $\partial_t (u,w) + \partial_x \tilde f(u,w) = \tilde
G(u,w)$ generates the operator $\tilde F$. The Cauchy problem
\begin{displaymath}
  \left\{
    \begin{array}{l}
      \partial_t u + \partial_x f(u) = G(t,u)
      \\
      u(0,x) = u_o(x)
    \end{array}
  \right.
\end{displaymath}
is solved by $ t \mapsto F_t(u_o,0)$ where $F_t(u_o,0)$ is given by
the first $n$ component of $\tilde F_{t} (u_o,0)$.

Again, the integral estimates~\ref{it:thm:main} and~\ref{it:thm2:main}
in Theorem~\ref{thm:main} ensure that the solutions constructed here
coincide with those in~\cite{CrastaPiccoli}.

\Section{Technical Proofs}
\label{sec:Tech}

\begin{proofof}{Proposition~\ref{prop:G}}
  The Lipschitz property is immediate. To prove the bound on the total
  variation, call $\Lip(h)$ the Lipschitz constant of $h$. Then, it is
  sufficient to compute:
  \begin{eqnarray*}
    \tv \left( Q*h(u) \right)
    & \leq &
    \sup \sum_i \int_\reali \norma{Q(y)} \, 
    \norma{h\left(u(x_i-y)\right) - h\left( u(x_{i-1}-y)\right)} \, dy
    \\
    & \leq &
    \Lip(h) \, \sup \sum_i \int_\reali \norma{Q(y)} \, 
    \norma{u(x_i-y) - u(x_{i-1}-y)} \, dy
    \\
    & \leq &
    \Lip(h) \, \sup \int_\reali \norma{Q(y)} \, 
    \sum_i \norma{u(x_i-y) - u(x_{i-1}-y)} \, dy
    \\
    & \leq &
    \Lip(h) \, \sup \int_\reali \norma{Q(y)} \, 
    \tv (u) \, dy
    \\
    & \leq &
    \Lip(h) \, \norma{Q}_{\L1} \, \tv (u) 
  \end{eqnarray*}
\end{proofof}

\subsection{Convective part}

Let $\lambda_1(u), \lambda_2(u), \ldots, \lambda_n(u)$ be the $n$ real
distinct eigenvalues of $D\!f(u)$, indexed so that $\lambda_j <
\lambda_{j+1}$ for all $j$ and $u$. The $j$-th right eigenvector is
$r_j(u)$ and we assume that $\norma{r_j(0)} =1$.

Let $\sigma \mapsto R_j(\sigma)(u)$ and $\sigma \mapsto
S_j(\sigma)(u)$ be respectively the rarefaction and the shock curve
exiting $u$. If the $j$-th field is linearly degenerate, then the
parameter $\sigma$ above is the arc-length. In the genuinely nonlinear
case, see~\cite[Definition~5.2]{BressanLectureNotes}, we choose
$\sigma$ so that for a suitable constant $k_j > 0$
\begin{displaymath}
  \begin{array}{rcl@{\quad\mbox{ and }\quad}rcl}
    \displaystyle
    \frac{\partial ~}{\partial\sigma}
    \lambda_j \left(R_j (\sigma)(u) \right) 
    & = &
    k_j
    &
    \displaystyle
    \frac{\partial R_j}{\partial \sigma} (0)(0) 
    & = &
    r_j(0)
    \\
    \\
    \displaystyle
    \frac{\partial ~}{\partial\sigma}
    \lambda_j \left(S_j (\sigma)(u) \right) 
    & = &
    k_j,
    &
    \displaystyle
    \frac{\partial S_j}{\partial \sigma} (0)(0)
    & = &
    r_j(0) \,.
  \end{array}
\end{displaymath}
The above choices were introduced in~\cite[\S~2]{AmadoriGuerra2002},
see also~\cite{BressanLectureNotes,BressanYangLiu}.

Introduce the $j$-Lax curve
\begin{displaymath}
  \sigma \mapsto \psi_j (\sigma) (u) =
  \left\{
    \begin{array}{c@{\qquad\mbox{ if }\quad}rcl}
      R_j(\sigma)(u) & \sigma & \geq & 0
      \\
      S_j(\sigma)(u) & \sigma & < & 0
    \end{array}
  \right.
\end{displaymath}
and define the map
\begin{displaymath}
  \Psi(\mathbf{\sigma})(u^-)
  =
  \psi_n(\sigma_n)\circ\ldots\circ\psi_1(\sigma_1)(u^-) \,.
\end{displaymath}
By~\cite[\S~5.3]{BressanLectureNotes}, given any two states $u^-,u^+
\in \Omega$ sufficiently close to $0$, there exists a vector
$(\sigma_1, \ldots, \sigma_n) = E(u^-,u^+)$ such that $u^+ =
\Psi(\mathbf{\sigma})(u^-)$.

Similarly, let ${\bf S}$ be defined by
\begin{displaymath}
  u^+
  =
  {\bf S}(\mathbf{\sigma})(u^-) =
  S_n(\sigma_n) \circ \ldots \circ S_1(\sigma_1) (u^-)
\end{displaymath}
as the gluing of the Rankine - Hugoniot curves.

For a sufficiently small $\delta_o$, let $u \in
\mathcal{U}_{\delta_o}$ be piecewise constant with finitely may jumps
sited in a finite set of points denoted by $\mathcal{I}(u)$.  Let
$\sigma_{x,i}$ be the strength of the $i$-th wave in the solution of
the Riemann problem for~\Ref{eq:HCL} with data $u(x-)$ and $u(x+)$.
i.e.~$(\sigma_{x,1}, \ldots, \sigma_{x,n}) = E\left( u(x-), u(x +)
\right)$. Obviously if $x\not\in \mathcal{I}(u)$ then
$\sigma_{x,i}=0$, for all $i=1,\ldots,n$.  As
in~\cite[\S~7.7]{BressanLectureNotes}, $\mathcal{A}(u)$ denotes the
set of approaching waves in $u$:
\begin{displaymath}
  \mathcal{A} (u) 
  = \left\{
    \begin{array}{c}
      \left(
        (x,i),(y,j)\right) \in \left( \mathcal{I}(u) \times \{1,\ldots,n\} 
      \right)^2
      \colon \\
      x < y \mbox{ and either } i > j \mbox{ or } i = j, 
      \mbox{ the $i$-th field}\\
      \mbox{is genuinely non linear, }
      \min \left\{ \sigma_{x,i}, \sigma_{y,j} \right\} <0\!  
    \end{array}
  \right\}
\end{displaymath}
while the linear and the interaction potential are
\begin{displaymath}
  V(u)
  =
  \sum_{x\in I(u)} \sum_{i=1}^n
  \modulo{\sigma_{x,i}}
  \,,\qquad
  Q(u)
  =
  \sum_{\left((x,i),(y,j)\right) \in \mathcal{A}(u)}
  \modulo{\sigma_{x,i}\sigma_{y,j}} \,.
\end{displaymath}
Moreover, let
\begin{equation}
  \Upsilon(u) = V(u) + C_0 Q (u)
  \label{def:ups}
\end{equation}
where $C_0>0$ is the constant appearing in the functional of the
wave--front tracking algorithm,
see~\cite[Proposition~7.1]{BressanLectureNotes}. Finally we define
\begin{equation}
  \overline{\mathcal{D}}_\delta
  =
  \mathrm{cl}
  \left\{
    u \in \L1\left(\reali,\reali^n\right) \colon
    u \hbox { is piecewise constant and } \Upsilon(u) < \delta
  \right\}
  \label{def:2.6}
\end{equation}
where the closure is in the strong $\L1$--topology. We remark for
later use that there exists a positive constant $c = c(\delta_o)$ with
$c \in \left]0, 1 \right[$, such that for all $\delta \in \left[0,
  \delta_o\right]$
\begin{displaymath}
  \mathcal{U}_\delta 
  \supseteq \overline{\mathcal{D}}_{c\delta}
  \supseteq \mathcal{U}_{c^2 \delta}
  \,.
\end{displaymath}

\begin{proposition}
  \label{prop:SRS}
  Let $f$ satisfy~\textbf{(F)}. Then, there exists a positive and
  suitably small $\bar\delta_o$ such that~\Ref{eq:HCL} generates a
  Standard Riemann Semigroup (SRS), with Lipschitz constant $L$,
  defined on the domain $\overline{\mathcal{D}}_\delta$, for all
  $\delta \in \left]0, \bar\delta_o \right[$.
\end{proposition}

We refer to~\cite[Chapters~7 and~8]{BressanLectureNotes} for the proof
of the above result as well as for the definition and further
properties of the SRS.

\begin{lemma}
  \label{lem:Taylor}
  Let $f$ satisfy~\textbf{(F)}, $\Omega$ be a sufficiently small
  neighborhood of the origin; $a,b \in \reali^n$ and $s \in \left[0,
    +\infty\right[$ be sufficiently small. Choose $u^-,v^- \in \Omega$
  and define
  \begin{displaymath}
    u^+ = u^- + s a \,,\ v^+ = v^- + s b \,.
  \end{displaymath}
  Then, there exist $\mathbf{\sigma^-}$, $\mathbf{\sigma^+}$ such that
  $v^- = \Psi (\mathbf{\sigma^-}) (u^-)$ and $v^+ = \Psi
  (\mathbf{\sigma^+}) (u^+)$. Moreover,
  \begin{equation}\label{SizeEstimate}
    \sum_{i=1}^n \modulo{\sigma^+_i - \sigma^-_i}
    \leq
    \O \cdot
    \left( \norma{a-b} + \sum_{i=1}^n \modulo{\sigma_i^-} \right)
    \cdot s\,.
  \end{equation}
  An entirely analogous result holds with $\Psi$ replaced by
  $\mathbf{S}$, i.e.~$v^- = \mathbf{S} (\mathbf{\sigma^-}) (u^-)$ and
  $v^+ = \mathbf{S} (\mathbf{\sigma^+}) (u^+)$.
\end{lemma}

\noindent The proof is an extension
of~\cite[Lemma~2.1]{AmadoriGuerra2002} and, hence, omitted.

\subsection{Source part}

Concerning the source term we have the following results.

\begin{proposition}
  \label{Lem:Source} Let $G$ satisfy~\textbf{(G)}. Then, for any
  $\delta \in \left]0,\delta_o\right[$ and $u\in\mathcal{U}_\delta$,
  the Cauchy Problem~\Ref{eq:Source} with initial data $u$ admits a
  solution $\Sigma_t u$ defined for $t \in [0,\widetilde T]$ where
  $\widetilde T =
  \min\left\{\frac{\delta_o-\delta}{\delta_oL_2+L_3},\frac{1}{L_1+1}
  \right\}$. Moreover the trajectory $\Sigma_tu$ has the following
  properties for all $t \in [0,\widetilde T]$ and
  $u,v\in\mathcal{U}_\delta$:
  \begin{eqnarray*}
    \displaystyle
    \norma{\Sigma_tu - \Sigma_tw}_{\L1} 
    & \leq &
    \displaystyle
    \norma{u-w}_{\L1} \, e^{L_1 t} \,
    \\
    \displaystyle
    \tv \left( \Sigma_t u \right)
    & \leq & 
    \displaystyle
    \delta + \left(\delta_o \ L_2+L_3\right)t.
  \end{eqnarray*}
\end{proposition}

The existence of $\Sigma_tu$ is a standard application of Banach Fixed
Point Theorem and the estimates follow from Gr\"onwall Lemma.

Since the computations on the convective part is mainly done on
piecewise constant approximate solutions, we need to approximate the
source term with piecewise constant functions.

Let $\PC(\reali;\reali^n)$ be the set of piecewise constant functions
in $\L1$. For any $N \in \naturali$, define the operator $\Pi_N \colon
\L1(\reali;\reali^n) \mapsto \PC(\reali;\reali^n)$ by
\begin{displaymath}
  \Pi_N (u)
  =
  N \sum_{k=-1-N^2}^{-1+N^2} \int_{k/N}^{(k+1)/N} u(\xi) \, d\xi \;
  \chi_{\strut \left]k/N, (k+1)/N \right]} \,.
\end{displaymath}

\begin{lemma}\label{lemmaconvergence}
  $\Pi_N$ is a linear operator with norm $1$. Moreover, $\tv \left(
    \Pi_N u \right) \leq 2 \tv(u)$ and for all $u \in
  \L1(\reali;\reali^n) \cap \BV(\reali;\reali^n)$, $\Pi_Nu \to u$ in
  $\L1 (\reali;\reali^n)$.
\end{lemma}

\begin{proof}
  Linearity and the estimate on the norm are immediate.
  \begin{eqnarray*}
    \tv ( \Pi_N u )
    & = &
    \norma{(\Pi_N u) (-N)}
    +
    \norma{(\Pi_N u) (N)}
    \\
    & &
    +\sum_{k=-N^2}^{N^2-1}
    \norma{(\Pi_N u) \left( (k+1)/N\right) - (\Pi_N u) (k/N)}
    \\
    & \leq &
    N \int_{-N-1/N}^{-N} \norma{u(\xi)}\, d\xi +
    N \int_{N - 1/N}^N \norma{u(\xi)}\, d\xi
    \\
    & &
    + N
    \sum_{k=-N^2}^{N^2-1}
    \norma{ \int_{k/N}^{(k+1)/N} u(\xi)\, d\xi -
      \int_{(k-1)/N}^{k/N} u(\xi)\, d\xi}
    \\
    & \leq &
    \tv \left( u; \left]-\infty, -N\right[\right)
    +
    \tv \left( u; \left]N - 1/N, +\infty \right[\right)
    \\
    & &
    + N \sum_{k=-N^2}^{N^2-1}
    \int_{k/N}^{(k+1)/N} \norma{u(\xi) - u(\xi-1/N)}\, d\xi
    \\
    & \leq &
    \tv \left( u; \left]-\infty, -N\right[\right)
    +
    \tv \left( u; \left]N-1/N, +\infty \right[\right)
    \\
    & &
    + \sum_{k=-N^2}^{N^2-1} \tv \left(u; \left](k-1)/N, (k+1)/N \right[ \right)
    \\
    & \leq &
    2 \, \tv (u) \,.
  \end{eqnarray*}
  Concerning the pointwise convergence $\Pi_N \to \Id$:
  \begin{eqnarray*}
    \norma{\Pi_Nu - u}_{\L1}
    & \leq &
    \int_{-\infty}^{-N-1/N} \norma{u(\xi)}\, d\xi
    +
    \int_{N}^{+\infty} \norma{u(\xi)} \, d\xi
    \\
    & &
    + N \sum_{k=-N^2-1}^{N^2-1}
    \int_{k/N}^{(k+1)/N} \int_{k/N}^{(k+1)/N}
    \norma{u(\xi) - u(x)} \, d\xi \, dx
    \\
    & \leq &
    \int_{-\infty}^{-N-1/N} \norma{u(\xi}\, d\xi
    +
    \int_{N}^{+\infty} \norma{u(\xi)} \, d\xi
    +
    \frac{1}{N}\, \tv(u)
  \end{eqnarray*}
  So that $\lim_{N\to+\infty} \norma{\Pi_Nu - u}_{\L1} = 0$.
\end{proof}

\begin{corollary}
  Let $G$ satisfy~\textbf{(G)}. Then, for any $N$, also $\Pi_N \circ
  G$ satisfies~\textbf{(G)} with the same Lipschitz constant $L_1$ and
  with $L_2,L_3$ replaced by $2L_2$ and $2L_3$.
\end{corollary}

To simplify the operator splitting algorithm, we substitute the
semigroup generated by~\Ref{eq:Source} with the Euler approximation
$P_t \colon \mathcal{U}_{\delta_o} \mapsto \L1 (\reali;\reali^n)$
defined by
\begin{displaymath}
  P_t u = u + t \; G(u) \,.
\end{displaymath}
It is immediate to prove that $P_t$ is $\L1$-Lipschitz with constant
$1+tL_1$. We often use below the estimate $1+L_1 t \leq e^{L_1 t}$. On
the other hand, observe that, despite the notation (useful in the
sequel), $P$ does not satisfy the semigroup composition law. Indeed we
can only say
\begin{displaymath}
  \norma{P_sP_tu-P_{s+t}u}_{\L1}
  \leq 
  \O \cdot s\cdot t \cdot \left( 1+ \norma{u}_{\L1} \right) \,.
\end{displaymath}

\subsection{Operator Splitting}

\begin{lemma}
  \label{Lemma2.6}
  Let $G$ attain values in $\PC(\reali;\reali^n)$. Then, there exist
  positive $\bar\delta_o$ and $T_o$ such that for all $s \in [0,T_o]$
  and for all piecewise constant $u \in \mathcal{D}_{\bar\delta_o}$
  \begin{eqnarray}
    \label{Vestimate}
    V \left( P_s u \right)
    & \leq &
    V(u) + \O  s  \left[ L_3+ V(u) \right]
    \\\label{Qestimate}
    Q \left( P_s u \right)
    & \leq &
    Q(u) + \O  s  \left[ L_3+ V^2(u)\right]
    \\\label{UPSestimate}
    \Upsilon \left( P_s u \right)
    & \leq &
    \Upsilon(u) + \O  s  \left[  L_3+ V(u) \right]
  \end{eqnarray}
\end{lemma}

\begin{proof}
  Let $u' = P_su$, so that $u' = u + s \, G(u)$. Call $\sigma'_{x,i}$
  (resp.~$\sigma_{x,i}$) the size of the $i$-wave in $u'$ (resp.~$u$)
  at the point $x$, observe that $\sigma'_{x,i}$
  (resp.~$\sigma_{x,i}$) vanishes whenever $x \not \in
  \mathcal{I}(u')$ (resp.~$x \not \in \mathcal{I}(u)$).  Then, by
  Lemma~\ref{lem:Taylor} and~\textbf{(G)}
  \begin{eqnarray*}
    V(u') - V(u)
    & \leq &
    \sum_{x \in \mathcal{I}(u') \cup \mathcal{I}(u)}
    \sum_{i=1}^n \modulo{\sigma_{x,i} - \sigma_{x,i}'}
    \\
    & \leq &
    \O \cdot s \cdot \sum_{x \in \mathcal{I}(u') \cup \mathcal{I}(u)}
    \left(
      \sum_{i=1}^n \modulo{\sigma_{x,i}} +
      \norma{\Delta \left( G(u)\right)(x)}
    \right)
    \\
    & \leq &
    \O \cdot s \cdot \left( V(u) + \tv \left(G(u)\right)\right)
    \\
    & \leq &
    \O \cdot s \cdot \left( L_3 + V(u)\right) \,.
  \end{eqnarray*}
  To derive~\Ref{Qestimate}, we observe that if $\left((x,i), (y,j)
  \right) \in \mathcal{A}(u')\backslash\mathcal{A}(u)$, then either
  $\sigma_{x,i}'\sigma_{x,i}\le 0$ or $\sigma_{y,j}'\sigma_{y,j}\le
  0$. Suppose, for instance, that $\sigma_{x,i}'\sigma_{x,i}\le 0$
  (the other case being similar) then by~\Ref{SizeEstimate} we obtain
  \begin{displaymath}
    \modulo{\sigma_{x,i}'} + \modulo{\sigma_{x,i}}
    = 
    \modulo{\sigma_{x,i}' - \sigma_{x,i}}
    \leq
    \O s \,
    \left( \sum_{j=1}^n \modulo{\sigma_{x,j}} + \norma{\Delta G(u)(x)} \right)
  \end{displaymath}
  and therefore
  \begin{equation}
    \label{eq:2.9}
    \sum_{\left((x,i),(y,j)\right)\in\mathcal{A}(u')\backslash\mathcal{A}(u)}
    \modulo{\sigma_{x,i} \sigma_{y,j}}
    \leq 
    \O s \left( V(u) + L_3\right) V(u) \,.
  \end{equation}
  Applying~\Ref{eq:2.9}, \Ref{Vestimate} and since $s,V(u),V(u')<<1$
  we finally get
  \begin{eqnarray*}
    Q(u') & = &
    \sum_{\mathcal{A}(u')} \modulo{\sigma'_{x,i} \sigma'_{y,j}}
    \\
    &\leq&
    \sum_{\mathcal{A}(u')} \left( \modulo{\sigma'_{x,i} -
        \sigma_{x,i}} \modulo{\sigma'_{y,j}} + \modulo{\sigma_{x,i}}
      \modulo{\sigma'_{y,j} -\sigma_{y,j}} \right)+
    \sum_{\mathcal{A}(u')} \modulo{\sigma_{x,i} \sigma_{y,j}}
    \\
    &\leq &
    \O
    s\left\{\left[V(u)+L_3\right]V(u')+V(u)\left[V(u)+L_3\right]
    \right\}
    \\
    & &
    +\sum_{\mathcal{A}(u')\setminus \mathcal{A}(u)}
    \modulo{\sigma_{x,i} \sigma_{y,j}} +
    \sum_{\mathcal{A}(u)} \modulo{\sigma_{x,i} \sigma_{y,j}}
    \\
    &\leq &
    \O \cdot s \cdot \left( V(u)^2 + L_3 \right) + Q(u) \,.
  \end{eqnarray*}
  Since $L_3$ is a possibly null constant which depends only on the
  system, we can say $L_3+V(u)=\O$. Finally, the latter estimate
  follows combining the previous results.
\end{proof}

\begin{corollary}\label{cor2.7}
  Let $\delta \in \left]0, \bar\delta_o\right[$ and assume that $G$
  satisfies~\textbf{(G)}. Then we have
  \begin{displaymath}
    P_s \overline{\mathcal{D}}_\delta
    \subseteq
    \overline{\mathcal{D}}_{\delta+\O s \left(2L_3+\delta\right)}
    \quad \mbox{ and } \quad
    P_s \mathcal{U}_{\delta}
    \subseteq
    \mathcal{U}_{\delta+s(L_2\delta+L_3)} \,.
  \end{displaymath}
  In particular take a constant $C\ge \O \left(2L_3 + \bar
    \delta_o\right)$, a number $\delta \in \left]0, \bar\delta_o
  \right[$ and time $T \in \left]0, T_o\right]$ such that $\delta +
  C\,T \leq \bar\delta_o$ then for any $t\in[0,T]$, $s\in[0,T-t]$ we
  have
  \begin{equation}\label{Deltaest}
    P_s \overline{\mathcal{D}}_{\delta+tC} \subseteq
    \overline{\mathcal{D}}_{\delta+C(t+s)}
  \end{equation}
\end{corollary}

\begin{proof}
  Fix $u \in \overline{\mathcal{D}}_\delta$ and an approximating
  sequence of piecewise constant function $u_k$ with $\Upsilon(u_k) <
  \delta$. The previous Lemma shows that
  \begin{eqnarray*}
    \Upsilon\left( u_k + s (\Pi_N\circ G) (u_k) \right)
    & < &
    \Upsilon(u_k) + \O s\left[2L_3+V(u_k)\right]
    \\
    & < &
    \delta+\O s\left[2L_3+\delta\right] \,.
  \end{eqnarray*}
  But $u_k+s\Pi_N\circ G(u_k)$ converges to $P_su$ as $k,N\rightarrow
  +\infty$ and so $P_s u \in \overline{\mathcal{D}}_{\delta+\O s
    \left(2L_3+\delta\right)}$. The proofs of the other inclusions are
  straightforward.
\end{proof}

Corollary~\ref{cor2.7} allows to define the domains appearing in
Theorem~\ref{thm:main} as
\begin{displaymath}
  \mathcal{D}_t = \overline{\mathcal{D}}_{\delta + Ct}
  \qquad
  \forall\, t \in [0, T] \,.
\end{displaymath}

Let $h \in \naturali$ and define
\begin{equation}
  \label{eq:DefP}
  F^s_{t} u =  S_{t-hs} (P_s \circ S_s)^h\, u
  \qquad
  t \in \left[hs, (h+1)s \right[ \,.
\end{equation}
In other words, in any interval $\left]hs, (h+1) s \right[$, we apply
the semigroup $S$.  In turn, at the times $t=hs$, $P_s$ is applied.

If a time $T$ and a $\delta \in \left]0, \bar\delta_o \right[$ are
chosen as in Corollary~\ref{cor2.7}, then $F^s_tu$ is defined up to
the time $T-\bar t$ for any $u \in \mathcal{D}_{\bar t}$ and $\bar t
\in [0,T]$.

Observe that for $t',t'' \in [0,T]$ with $t' + t'' \in [0,T]$, the
following inclusion holds:
\begin{equation}
  \label{eq:inclusion}
  F^s_{t'} \mathcal{D}_{t''} \subseteq \mathcal{D}_{t'+t''} \,.
\end{equation}
Note that $F^s_t$ is $\L1$-Lipschitz with Lipschitz constant bounded
by $L\cdot L^{t/s} e^{L_1t}$, with $L$ as in
Proposition~\ref{prop:SRS} and $L_1$ as in~\textbf{(G)}. The following
theorem shows that the Lipschitz constant of $F^s_t$ actually is
bounded from above by a quantity independent from $s$.

\begin{theorem}
  \label{thm:Lipschitz}
  Let $f$ satisfy~\textbf{(F)} and $G$ satisfy~\textbf{(G)}.  If
  $\bar\delta_o$ is chosen sufficiently small and $\delta,T$ are
  chosen as in Corollary~\ref{cor2.7}, then there exists a constant
  $\mathcal{L}$ such that for all $\bar t \in [0,T]$, $u, w \in
  \mathcal{D}_{\bar t}$ and $t \in [0, T-\bar t]$, we have
  \begin{equation}\label{lip1propo}
    \norma{F^s_t u - F^s_t w}_{\L1}
    \leq
    \mathcal{L} \cdot \norma{u-w}_{\L1}\,.
  \end{equation}
\end{theorem}

\begin{proof}
  The proof can be carried out following essentially the same line
  used in~\cite[Theorem~4.2]{AmadoriGuerra2002}. Therefore we only
  outline it. The key ingredient is the Liu \& Yang functional which,
  unfortunately, is defined only on piecewise constant
  $\varepsilon$--approximations of the trajectories of the convective
  part (see~\cite{BressanLectureNotes}). Therefore we first suppose
  that the source $G(u)$ attains values in
  $\PC\left(\reali,\reali^n\right)$. For this kind of sources one can
  prove~\Ref{lip1propo} following almost exactly the proof
  of~\cite[Theorem~4.2]{AmadoriGuerra2002}. The only difference is the
  {\it global} nature of our source term. But the non locality of
  $G(u)$ can easily be tackled using Lemma~\ref{lem:Taylor},
  integrating and summing up the pointwise estimates obtained
  in~\cite[Theorem~4.2]{AmadoriGuerra2002} and finally using {\bf (G)}
  similarly to what we have done in a detailed way in the proof of
  Lemma~\ref{Lemma2.6}.

  Once we obtained~\Ref{lip1propo} for piecewise constant source
  terms, we apply it to the source term $\Pi_N\circ G(u)$. If we
  denote by $F^{s,N}_tu$ the trajectory defined by~\Ref{eq:DefP} with
  $\Pi_N\circ G(u)$ in place of $G(u)$, it is quite easy to see that
  Lemma~\ref{lemmaconvergence} implies the strong convergence of
  $F^{s,N}_tu$ to $F^{s}_tu$ proving the theorem for a general source
  $G$.
\end{proof}

We need below the following estimates concerning the dependence of
$F^s$ on time.

\begin{lemma}
  Let $T$ and $\delta$ be as in corollary~\ref{cor2.7}. Then, for all
  $u \in \mathcal{D}_{\bar t}$ and all $t \in [0, T-\bar t]$,
  \begin{eqnarray}
    \label{eq:BoundF}
    \norma{F^s_t u}_{\L1}
    & \leq &
    \O \cdot \left( 1+ \norma{u}_{\L1} \right)
    \\
    \label{eq:BoundP}
    \norma{(P_s)^k u}_{\L1}
    & \leq &
    \O \cdot \left( 1+ \norma{u}_{\L1} \right)
    \\
    \label{eq:LipF}
    \norma{F^s_t u - u}_{\L1}
    & \leq &
    \O \cdot t \cdot \left( 1+ \norma{u}_{\L1} \right)
    \\
    \label{eq:LipP}
    \norma{(P_s)^k u - u}_{\L1}
    & \leq &
    \O \cdot k \cdot s \cdot \left( 1+ \norma{u}_{\L1} \right) \,.
  \end{eqnarray}
\end{lemma}

\begin{proof}
  Consider first~\Ref{eq:BoundP} with $k=1$. By~\textbf{(G)} and the
  properties of $S$
  \begin{eqnarray*}
    \norma{P_s u}_{\L1}
    & = &
    \norma{u + s \, G(u) - s\, G(0) + s\,G(0)}_{\L1}
    \\
    & \leq &
    \norma{u}_{\L1} + s \, L_1 \norma{u}_{\L1} + s \, \norma{G(0)}_{\L1}
    \\
    & \leq & e^{sL_1} \left( \tilde c \, s + \norma{u}_{\L1} \right)
    \\
    \norma{S_s u}_{\L1}
    & \leq &
    \norma{S_s u -u}_{\L1} + \norma{u}_{\L1}
    \\
    & \leq &
    L\, s + \norma{u}_{\L1}
    \\
    & \leq &
    e^{Ls} \left( \tilde c \, s + \norma{u}_{\L1} \right)
  \end{eqnarray*}
  where $\tilde c = \max \{L, \norma{G(0)}_{\L1}\}$ and $L$ is as in
  Proposition~\ref{prop:SRS}. Proceed now by induction on $k$ and, for
  $k\,s \in \left[0, T-\bar t \right[$,
  \begin{eqnarray*}
    \norma{F^s_{ks} u}_{\L1}
    & \leq &
    e^{(L+L_1) k s} \left( 2\,\tilde c \, k\,s + \norma{u}_{\L1} \right)
    \\
    & \leq &
    e^{(L+L_1) k s} \left( 2\,\tilde c \, T+ \norma{u}_{\L1} \right)
    \\
    & \leq &
    \O \cdot \left( 1 + \norma{u}_\L1 \right) \,.
  \end{eqnarray*}
  Therefore, for $t \in \left[0, T - \bar t\right[$ and $ \bar k =
  \left[ \frac{t}{s}\right]$, we have
  \begin{eqnarray*}
    \norma{F^s_t u}_{\L1}
    & = &
    \norma{S_{t-ks} F^s_{ks} u}_{\L1}
    \\
    & \leq &
    L\, s + \O \cdot \left( 1 + \norma{u}_{\L1} \right)
    \\
    & \leq &
    \O \cdot \left( 1 + \norma{u}_{\L1} \right) \,.
  \end{eqnarray*}
  The estimate~\Ref{eq:BoundP} is obtained similarly. Passing to the
  Lipschitz estimates, for $(k+1)s \in \left[0, T - \bar t \right[$,
  by~\Ref{eq:BoundF}
  \begin{eqnarray*}
    \norma{F^s_{(k+1)s} u - u}_{\L1}
    & = &
    \norma{F^s_{s} \, F^s_{ks} u - F^s_{ks}u}_{\L1} +
    \norma{F^s_{ks}u - u}_{\L1}
    \\
    & \leq &
    \norma{S_s F^s_{ks} u + s G(S_s F^s_{ks} u) - F^s_{ks}u}_{\L1} +
    \norma{F^s_{ks}u - u}_{\L1}
    \\
    & \leq &
    s \cdot \left( 
      L + L_1 \, \norma{S_s F^s_{ks} u}_{\L1} + \norma{G(0)}_{\L1} 
    \right) +
    \norma{F^s_{ks}u - u}_{\L1}
    \\
    & \leq &
    \O \cdot s \cdot \left( 1 + \norma{u}_{\L1} \right)  +
    \norma{F^s_{ks}u - u}_{\L1} \,.
  \end{eqnarray*}
  By induction, $\norma{F^s_{ks} u -u}_{\L1} \leq \O \cdot k \cdot s
  \cdot \left( 1+\norma{u}_{\L1} \right)$.

  We are left with the case $\frac{t}{s} \not\in \naturali$. If $t \in
  \left[0, s \right[$ we have
  \begin{displaymath}
    \norma{F^s_t u - u}_{\L1} = \norma{S_t u -u}_{\L1} \leq \O \cdot t \,,
  \end{displaymath}
  while if $t \geq s$, so that $\bar k= \left[ \frac{t}{s}\right] \geq
  1$,
  \begin{eqnarray*}
    \norma{F^s_t u - u}_{\L1}
    & \leq &
    L\, s + \norma{F^s_{\bar ks} u - u}_{\L1}
    \\
    & \leq &
    L \, s + \O \cdot \bar k \cdot s \cdot \left(1 + \norma{u}_{\L1} \right)
    \\
    & \leq &
    \O \cdot t \cdot \left(1 + \norma{u}_{\L1} \right) \,.
  \end{eqnarray*}
  The proof of~\Ref{eq:LipP} is entirely similar.
\end{proof}

The next step consists in showing the convergence of $F^s$ as $s$
tends to zero. This result will be obtained with the help of the
commutation relation~\Ref{commutation} which we will show to be true
for $S_t$ and $P_t$. We need the following result that is an easy
consequence of~\cite[Remark~4.1]{AmadoriGuerra2002}.

\begin{proposition}\label{proprem}
  Take $u,v,\omega\in\mathcal{D}_{T}$, then, for any $t_1<t_2$ one has
  the estimate
  \begin{eqnarray}
    \label{eq:2.15}
    \norma{S_{t_2} w - S_{t_2} u - \omega}_{\L1}
    & \leq &
    L \cdot \norma{S_{t_1} w - S_{t_1} u - \omega}_{\L1}
    \\
    \nonumber
    & &
    + \O \cdot (t_2-t_1) \cdot \tv(\omega) \,.
  \end{eqnarray}
\end{proposition}

Proposition~\ref{proprem} implies the following commutation result:

\begin{theorem}
  \label{thm:2.10}
  Let $\delta$ and $T$ be in Corollary~\ref{cor2.7}. For any $\bar t
  \in \left[0,T \right[$, $u \in \mathcal{D}_{\bar t}$ and $t \in
  \left[0,T-\bar t\right]$, we have the estimate
  \begin{equation}
    \label{commutation2}
    \norma{S_{t}P_t u-P_tS_{t} u}_{\L1} \leq \O \cdot t^2 \, .
  \end{equation}
\end{theorem}

\begin{proof}
  Since $u + t G(u)$, $tG(S_tu)$ and $u$ all belong to
  $\mathcal{D}_T$, we can apply~\Ref{eq:2.15} with $w = u + t\, G(u)$
  and $\omega = t \, G \left(S_t u\right)$ to obtain:
  \begin{eqnarray*}
    \norma{S_{t}P_t u-P_tS_{t} u}_{\L1}
    & = &
    \norma{
      S_{t}\left[u + t \, G(u)\right]- 
      S_{t} u - t\, G\left(S_t u\right)}_{\L1}
    \\
    &\leq &
    L \cdot \norma{S_{0}\left[u + t\, G(u)\right]
      -
      S_{0} u - t \, G \left(S_t u\right)}_{\L1}
    \\
    & &
    +
    \O \cdot t \cdot \tv\left[t \, G\left(S_t u\right)\right]
    \\
    &\leq &
    L \cdot t \cdot \norma{G(u) - G\left(S_t u\right)}_{\L1} + 
    \O \cdot t^2
    \\
    &\leq &
    L^2 \cdot L_1 \cdot t^2 + \O \cdot t^2 \,.
  \end{eqnarray*}
\end{proof}

Now we show that~\Ref{commutation2} and the uniform Lipschitz
property~\Ref{lip1propo} of the approximations imply the existence of
a ``tangent vector'' and the strong convergence of the approximations.
We will show that these two conditions are enough and that there is no
need to use again the almost decreasing functional as was done
in~\cite[Lemma 5.1]{AmadoriGuerra2002}.

\begin{proposition}
  Let $\delta$ and $T$ be as in Corollary~\ref{cor2.7}. For any $\bar
  t \in \left[0,T \right[$, $u \in \mathcal{D}_{\bar t}$, $t \in
  \left[0,T-\bar t\right]$ and $s,s' \in \left]0, t^2\right]$, we have
  \begin{equation}
    \label{firstof2.16}
    \begin{array}{rcl}
      \displaystyle
      \norma{F^s_t u-S_tP_{t} u}_{\L1}^{\phantom{\vert}}
      & \leq &
      \displaystyle
      \O \cdot \left( 1 + \norma{u}_{\L1} \right)\cdot t^2
      \\
      \displaystyle
      \norma{F^s_t u-P_tS_{t} u}_{\L1}^{\phantom{\vert}}
      & \leq &
      \displaystyle
      \O \cdot \left( 1 + \norma{u}_{\L1} \right)\cdot t^2
      \\
      \displaystyle
      \norma{F^s_t u- F^{s'}_tu}_{\L1}^{\phantom{\vert}}
      & \leq &
      \displaystyle
      \O  \cdot \left( 1 + \norma{u}_{\L1} \right)\cdot t^2 \,.
    \end{array}
  \end{equation}
\end{proposition}

\begin{proof}
  We prove only the first inequality in~\Ref{firstof2.16}, the other
  two inequalities being consequences of this one and of
  Theorem~\ref{thm:2.10}. For any integer $k \in [1, T/s]$ define
  \begin{displaymath}
    \rho_k(s)=\sup_{u \in \mathcal{D}_{T-ks}}
    \norma{(P_s)^kS_s u - S_s(P_s)^ku}_{\L1}.
  \end{displaymath}
  Now, for $u\in \mathcal{D}_{T-(k+1)s}$ (and hence $P_s u \in
  \mathcal{D}_{T-ks}$), we can compute
  \begin{eqnarray*}
    & &
    \norma{(P_s)^{k+1} S_su-S_s(P_s)^{k+1} u}_{\L1} \leq
    \\
    & \leq &
    \norma{(P_s)^kP_sS_su-(P_s)^kS_sP_su}_{\L1}
    +
    \norma{(P_s)^kS_sP_su-S_s(P_s)^kP_su}_{\L1}
    \\
    &\leq &
    e^{T L_1} \norma{P_sS_su-S_sP_su}_{\L1}
    +
    \sup_{w\in\mathcal{D}_{T-ks}}
    \norma{(P_s)^kS_sw-S_s(P_s)^kw}_{\L1}
    \\
    &\leq &
    e^{T L_1} \rho_1(s)+\rho_k(s) \,.
  \end{eqnarray*}
  And hence $\rho_{k+1} (s) \leq e^{TL_1} \rho_1(s) + \rho_k(s)$ that,
  by induction, gives
  \begin{displaymath}
    \rho_k(s)\le e^{TL_1} k \rho_1(s) \,.
  \end{displaymath}
  Now define
  \begin{displaymath}
    \bar\rho_k(s)
    =
    \sup_{u\in\mathcal{D}_{T-ks}}
    \norma{F^s_{ks}u-S_{ks}(P_s)^ku}_{\L1}.
  \end{displaymath}
  Again we can compute for $u \in \mathcal{D}_{T-(k+1)s}$
  \begin{eqnarray*}
    & &
    \norma{F^s_{(k+1)s} u-S_{(k+1)s}(P_{s})^{k+1} u}_{\L1} \leq 
    \\
    &\leq &
    \norma{F^s_{ks} P_s S_s u - F^s_{ks} S_s P_s u}_{\L1} 
    +
    \norma{F^s_{ks} S_s P_s u - S_{ks}(P_s)^k S_s P_s u}_{\L1}
    \\
    & &
    +
    \norma{S_{ks}(P_s)^kS_sP_su-S_{ks}S_s(P_s)^kP_su}_{\L1}
    \\
    &\leq &
    \mathcal{L}\rho_1(s)
    +
    \bar\rho_k(s)+L\norma{(P_s)^kS_sP_su-S_s(P_s)^kP_su}_{\L1}
    \\
    &\leq &
    \mathcal{L}\rho_1(s)+ \bar\rho_k(s)+L\rho_k(s)
    \\
    &\leq &
    \left(\mathcal{L}+Le^{L_1 T} k\right) \rho_1(s) + \bar\rho_k(s) \,.
  \end{eqnarray*}
  Therefore we have $\bar\rho_{k+1}(s) \leq \left(\mathcal{L} + L \,
    e^{L_1 T} k\right) \rho_1(s) + \bar\rho_k(s)$ which gives,
  together to $\rho_1(s) = \bar\rho_1(s)$, by induction
  \begin{displaymath}
    \bar\rho_k(s)
    \leq
    \left( \mathcal{L} + L \, e^{L_1T} k\right) \, k\, \rho_1(s).
  \end{displaymath}
  Fix now $t\in [0,T-\bar t\,]$, take $s \in \left]0, t^2\right[$ and
  define $\hat k = \left[\frac{t}{s}\right]$. We have for all $u \in
  \mathcal{D}_{\bar t} \subset \mathcal{D}_{T-\hat k s}$
  \begin{eqnarray*}
    \norma{F^s_{t} u-S_{t}(P_s)^{\hat k} u}_{\L1}
    & \leq &
    \O \cdot ( t - \hat k\, s) + \norma{F^s_{\hat k s} u - 
      S_{\hat k s} (P_s)^{\hat k} u}_{\L1}
    \\
    & \leq &
    \O \cdot (t - \hat k s)
    + 
    \left( \mathcal{L} + L \, e^{L_1 T} \frac{t}{s} \right)
    \frac{t}{s} \, \rho_1(s)
    \\
    & \leq &
    \O \cdot s + \O \cdot \left(t\,s + t^2\right) \frac{\rho_1(s)}{s^2}
    \\
    & \leq &
    \O \cdot t^2 \,,
  \end{eqnarray*}
  where the last inequality is a consequence to the fact that
  $\frac{\rho_1(s)}{s^2}$ is bounded because of~\Ref{commutation2}.

  We are left to prove that
  \begin{displaymath}
    \norma{(P_s)^{\hat  k} u - P_t u}_{\L1} 
    \leq 
    \O \cdot \left( 1 + \norma{u}_{\L1} \right) \cdot t^2 \,.
  \end{displaymath}
  By~\Ref{eq:LipP},
  \begin{eqnarray*}
    \norma{(P_s)^{k+1}u - P_{(k+1)s}u}_{\L1}
    & = &
    \norma{P_s (P_s)^k u - P_{k s}u - s \, G(u)}_{\L1}
    \\
    & = &
    \norma{(P_s)^k u - P_{k s}u + s \, G\left( (P_s)^k u\right) - s \, G(u)}_{\L1}
    \\
    & \leq &
    \norma{(P_s)^k u - P_{k s}u}_{\L1} + 
    \O \cdot \left( 1 + \norma{u}_{\L1} \right) \cdot k \, s^2
  \end{eqnarray*}
  We thus recursively obtain
  \begin{displaymath}
    \norma{(P_s)^k - P_{ks}u}_{\L1} 
    \leq 
    \O \cdot \left( 1 + \norma{u}_{\L1} \right) \cdot k^2 \, s^2 \,.
  \end{displaymath}
  Therefore
  \begin{eqnarray*}
    \norma{(P_s)^{\hat  k} u - P_t u}_{\L1}
    & \leq &
    \norma{(P_s)^{\hat  k} u - P_{\hat k s} u}_{\L1} + 
    \norma{P_{\hat k s} u - P_t u}_{\L1}
    \\
    & \leq &
    \O \cdot \left( 1 + \norma{u}_{\L1} \right) \cdot \hat k^2 \, s^2 
    + \O \cdot s
    \\
    & \leq &
    \O \cdot \left( 1 + \norma{u}_{\L1} \right) \cdot t^2 \,.
  \end{eqnarray*}
\end{proof}

Now we prove the convergence of the approximations and the
characterization of the tangent vector, i.e.~\ref{it:tg}. in
Theorem~\ref{thm:main}.

\begin{theorem}
  \label{thm:Conv}
  Let $\delta$ and $T$ be as in Corollary~\ref{cor2.7}. For any $\bar
  t \in \left[0,T\right[$, $u \in \mathcal{D}_{\bar t}$ and $t \in
  \left[0,T-\bar t\right]$ the sequence $F^s_tu$ converges in $\L1$ as
  $s\rightarrow 0$ to a limit trajectory $F_tu$ which satisfies the
  tangency conditions
  \begin{equation}
    \label{eq:2.25}
    \begin{array}{rcl}
      \displaystyle
      \norma{F_tu-S_tP_tu}_{\L1} 
      & \leq &
      \displaystyle
      \O \cdot \left( 1 + \norma{u}_{\L1} \right) \cdot t^2_{\phantom{\vert}}
      \\
      \displaystyle
      \norma{F_tu-P_tS_tu}_{\L1}^{\phantom{\vert}}
      & \leq &
      \displaystyle
      \O \cdot \left( 1 + \norma{u}_{\L1} \right) \cdot t^2 \,.
    \end{array}
  \end{equation}
\end{theorem}

\begin{proof}
  Because of~\Ref{firstof2.16}, we need only to show that
  $s\rightarrow F^s_tu$ is a Cauchy sequence in $\L1$ as $s\rightarrow
  0$. Fix $\varepsilon>0$ arbitrary. Then choose $0 = t_0 < t_1 <
  \ldots < t_{N-1} < t_N = t$ so that $t_{i}-t_{i-1} < \epsilon$ for
  $i=1,\ldots,N$. Then observe that Definition~\Ref{eq:DefP},
  Theorem~\ref{thm:Lipschitz} and~\Ref{eq:BoundF} imply that $F^s$
  satisfies an approximated semigroup condition:
  \begin{displaymath}
    \norma{F_{t_1}^s
      F_{t_2}^s u - F_{t_1+t_2}^s u}_{\L1}
    = 
    \O \cdot \left( 1 + \norma{u}_{\L1} \right) \cdot s \,.
  \end{displaymath}
  Therefore, for any $0<s,s'< \min_{i=1\ldots N}\left\{
    (t_{i}-t_{i-1})^2\right\}$, we can compute
  \begin{eqnarray*}
    \norma{F_t^{s'}u-F_t^su}_{\L1}
    & \leq &
    \sum_{i=1}^{N}
    \norma{F^s_{t-t_{i}}F_{t_i}^{s'}u-F_{t-t_{i-1}}^sF_{t_{i-1}}^{s'}u}_{\L1}
    \\
    & \leq &
    \sum_{i=1}^{N}
    \norma{F^s_{t-t_{i}}F_{t_i}^{s'}u-F_{t-t_{i}}^sF_{t_i-t_{i-1}}^sF_{t_{i-1}}^{s'}u}_{\L1}
    \\
    & &
    + \O\left( 1 + \norma{u}_{\L1} \right)
    Ns
    \\
    & \leq &
    \mathcal{L} \sum_{i=1}^{N}
    \norma{F_{t_i}^{s'}u-F_{t_i-t_{i-1}}^sF_{t_{i-1}}^{s'}u}_{\L1}
    +\O\left( 1 + \norma{u}_{\L1} \right)
    Ns
    \\
    &\leq &
    \mathcal{L} \sum_{i=1}^{N}
    \norma{F_{t_i-t_{i-1}}^{s'}F_{t_{i-1}}^{s'}u -
      F_{t_i-t_{i-1}}^sF_{t_{i-1}}^{s'}u}_{\L1}
    \\
    & &
    + \O \left( 1 + \norma{u}_{\L1} \right) \, N \, (s+s')
    \\
    & \leq &
    \O \left( 1 + \norma{u}_{\L1} \right) 
    \left( \sum_{i=1}^{N}(t_i-t_{i-1})^2 + N\, (s+s') \right)
    \\
    & \leq &
    \O \left( 1 + \norma{u}_{\L1} \right) 
    \left( \epsilon \, t + N \, (s+s') \right)\,.
  \end{eqnarray*}
  And, finally, as $s,s'\rightarrow 0$ we get
  \begin{displaymath}
    \limsup_{s,s'\rightarrow 0} \norma{F_t^{s'} u - F_t^s u}_{\L1}
    \leq
    \O \, \left( 1 + \norma{u}_{\L1} \right) \, \varepsilon \, t
  \end{displaymath}
  which proves the Theorem because of the arbitrariness of
  $\varepsilon$.
\end{proof}

The limit trajectory thus obtained satisfies~(\ref{it:thm:semi}) in
Theorem~\ref{thm:main}, as can be seen passing to the limit $s \to 0$
in~\Ref{eq:inclusion} and in the approximate semigroup condition
\begin{displaymath}
  \norma{F^s_{t_1} F^s_{t_2} u - F^s_{t_1+t_2} u}_{\L1}
  \leq \O \cdot \left( 1 + \norma{u}_{\L1} \right) \cdot s \,.
\end{displaymath}
Taking the same limit in~\Ref{lip1propo}, we prove the former
inequality in~\Ref{eq:error}. To prove the latter estimate, observe
that by~\Ref{eq:LipF}
\begin{displaymath}
  \norma{F_t u - u}_{\L1} 
  \leq 
  \O \cdot \left(1 + \norma{u}_{\L1} \right) \cdot t 
\end{displaymath}
while, for $t_2 > t_1$, the semigroup property implies
\begin{eqnarray*}
  \norma{F_{t_2} u - F_{t_1} u}_{\L1}
  & \leq &
  \norma{F_{t_2 - t_1} F_{t_1} u - F_{t_1} u}_{\L1}
  \\
  & \leq &
  \O \cdot \left( 1 + \norma{F_{t_1} u}_{\L1} \right) \cdot (t_2 - t_1) 
  \\
  & \leq &
  \O \cdot \left( 1 + \norma{u}_{\L1} \right) \cdot (t_2 - t_1) \,.
\end{eqnarray*}
Assertion~(\ref{it:tg}) follows from~\Ref{eq:2.25}.

We pass now to~(\ref{it:3}). The trajectory $t \mapsto F_t u$ is a
weak entropic solution of~\Ref{eq:Main}. This can be proved using the
properties of the approximate solutions constructed above, as
in~\cite{AmadoriGosseGuerra, AmadoriGuerra2002, CrastaPiccoli}. Here,
we prefer to exploit the tangent vector provided by
theorem~\ref{thm:Conv}

\begin{corollary}
  \label{cor:weak}
  Let $\delta$ and $T$ be as in Corollary~\ref{cor2.7}. For any $\bar
  t \in \left[0, T \right[$, $u \in \mathcal{D}_{\bar t}$ and $t \in
  [0,T-\bar t\,]$ the trajectory $t \mapsto F_t u$ is a weak entropic
  solution of~\Ref{eq:Main}.
\end{corollary}

\noindent For the definition of weak entropic solutions of a balance
law, refer to~\cite{DafermosBook}, \cite[(2.16)
and~(2.19)]{DafermosHsiao} or~\cite[\S~6]{CrastaPiccoli}.

\begin{proofof}{Corollary~\ref{cor:weak}}
  We show below only the entropy inequality, since the proof that $t
  \mapsto F_t u$ is a weak solution is entirely similar.

  Observe that, by~\Ref{eq:2.25}
  \begin{eqnarray*}
    F_t u 
    & = &
    S_t u + \O \cdot \left( 1 + \norma{u}_{\L1} \right) \cdot t
    \\
    F_t u 
    & = & 
    S_t u + t \, G(u) + \O \cdot \left( 1 + \norma{u}_{\L1} \right) \cdot t^2
  \end{eqnarray*}
  in $\L1$. Let $(\eta,q)$ be an entropy-entropy flux pair and $\phi
  \in \Cc1$ be a non negative test function.  Fix a positive
  $\epsilon$ and denote $I_i = \left[ i \epsilon, (i+1) \epsilon
  \right[ \times \reali$ for $i \in \naturali$. By the properties of
  $S$, $\partial_t \eta(S_t u) + \partial_x q(S_t u) \leq 0$ in the
  sense of distribution, using the Divergence Theorem we get
  \begin{eqnarray*}
    & &
    \int_0^T \!\! \int_\reali 
    \left( 
      \eta(F_t u) \partial_t \phi + q(F_t u) \partial_x \phi 
    \right) \, dx \, dt
    \quad =
    \\
    & = &
    \sum_i \int \!\!\! \int_{I_i} 
    \left( 
      \eta(F_t u) \partial_t \phi + q(F_t u) \partial_x \phi 
    \right) \, dx \, dt
    \\
    & = &
    \sum_i \left(
      \int \!\!\!\int_{I_i} 
      \left( 
        \eta(S_{t-i\epsilon} F_{i\epsilon} u) \partial_t \phi + 
        q(S_{t-i\epsilon} F_{i\epsilon} u) \partial_x \phi 
      \right) \, dx \, dt
    \right)
    \\
    & &
    + \O \left( 1 + \norma{u}_{\L1} \right) \epsilon
    \\
    & \geq &
    \sum_i \left(
      \int_{\reali} 
      \left( 
        \eta(S_{\epsilon} F_{i\epsilon} u) \, \phi\left( (i+1)\epsilon,x \right)
        - 
        \eta(F_{i\epsilon} u) \, \phi(i\epsilon,x)
      \right) \, dx
    \right)
    \\
    & &
    + \O \left( 1 + \norma{u}_{\L1} \right) \epsilon    \\
    & = &
    \sum_i \left(
      \int_{\reali} 
      \left( 
        \eta(S_{\epsilon} F_{i\epsilon} u) 
        - 
        \eta(F_{(i+1)\epsilon} u) 
      \right) \, \phi\left( (i+1)\epsilon,x \right)\, dx
    \right)
    \\
    & &
    + \O \left( 1 + \norma{u}_{\L1} \right) \epsilon    \\
    & = &
    \sum_i \left(
      \int_{\reali} 
      \left( 
        \eta(S_{\epsilon} F_{i\epsilon} u) 
        - 
        \eta(F_{\epsilon}F_{i\epsilon} u) 
      \right) \, \phi\left( (i+1)\epsilon,x \right)\, dx
    \right)
    \\
    & &
    + \O \left( 1 + \norma{u}_{\L1} \right) \epsilon    \\
    & = &
    \sum_i
    \int_{\reali} 
    \left( 
      \eta(S_{\epsilon} F_{i\epsilon} u) 
      - 
      \eta\left(
        S_{\epsilon}F_{i\epsilon} u + \epsilon G(F_{i\epsilon} u)
      \right) 
    \right) \, \phi\left( (i+1)\epsilon,x \right)\, dx
    \\
    & &
    + \O \left( 1 + \norma{u}_{\L1} \right) \epsilon
    \\
    & = &
    \sum_i
    \int_{\reali} 
    \left( 
      - D\eta( F_{i\epsilon} u) 
      \, \epsilon \, G( F_{i\epsilon} u)
      \, \phi\left( (i+1)\epsilon,x \right)\right)\, dx
    \\
    & &
    + \O \left( 1 + \norma{u}_{\L1} \right) \epsilon
    \\
    & = &
    \sum_i
    \int \!\!\!\int_{I_i} 
    \left( - D\eta( F_t u) \, G( F_t u) \, \phi( t,x) \right)\, dx \, dt
    \\
    & &
    + \O \left( 1 + \norma{u}_{\L1} \right) \epsilon    \\ 
    & = &
    \int_0^T \!\! \int_\reali 
    - D\eta( F_t u) \, G( F_t u) \, \phi( t,x) \, dx \, dt
    + \O \left( 1 + \norma{u}_{\L1} \right) \epsilon
  \end{eqnarray*}
  By the arbitrariness of $\epsilon$, we conclude with the
  distributional inequality
  \begin{displaymath}
    \partial_t \eta(F_tu) + \partial_x q(F_tu) - D\eta
    (F_t u) \, G(F_tu) \leq 0 \,.
  \end{displaymath}

\end{proofof}

Now we show a result on the dependence of the solution with respect to
the source term.

\begin{theorem}
  \label{thm:DipLip}
  Let $f_1,f_2$ satisfy~\textbf{(F)} and $G_1,G_2$
  satisfy~\textbf{(G)}. Call $F^1,F^2$ the corresponding semigroups
  and assume they are defined on a common family of domains
  $\mathcal{D}_{t}$. Then, for any $\bar t \in \left[0, T \right[$, $u
  \in \mathcal{D}_{\bar t}$ and $t \in \left[0 , T-\bar t\right]$, the
  Lipschitz estimate~\Ref{eq:Uffa} holds.

  Moreover, fix a flux $f$ satisfying~\textbf{(F)} and a sequence of
  source terms $G_k$ satisfying~\textbf{(G)}. If $G_k$ converges
  pointwise to $G$, then the corresponding semigroups $F^{k}$ converge
  pointwise to the semigroup $F$ generated by $G$.
\end{theorem}

\begin{proof}
  We apply the well known integral estimate,
  see~\cite[Theorem~2.9]{BressanLectureNotes}:
  \begin{displaymath}
    \norma{F^{1}_t u - F^{2}_t u}_{\L1}
    \leq
    \mathcal{L}
    \int_0^t \liminf_{\theta\rightarrow 0}
    \frac{\norma{F^{2}_{\theta+\tau}u - F^{1}_{\theta}F^{2}_{\tau}u}_{\L1}}{\theta} 
    \; d\tau
  \end{displaymath}
  By~\Ref{eq:2.25} and using the stability
  result~\cite[Corollary~2.5]{BianchiniColombo}, with obvious notation
  we have
  \begin{eqnarray*}
    & &
    \liminf_{\theta \to 0}
    \frac{1}{\theta}
    \norma{F^{2}_{\theta} F^{2}_{\tau}u-F^{1}_{\theta}F^{2}_{\tau}u}_{\L1}
    \, \leq
    \\
    & \leq &
    \liminf_{\theta \to 0}
    \frac{1}{\theta}
    \norma{S^2_\theta F^2_\tau u + \theta G_2(F^{2}_{\tau} u) 
      - S^1_\theta F^2_\tau u - \theta \, G_1(F^{2}_{\tau}u)}_{\L1}
    \\
    & \leq &
    \liminf_{\theta \to 0}
    \frac{1}{\theta}
    \norma{S^2_\theta F^2_\tau u - S^1_\theta F^2_\tau u }_{\L1}
    +
    \norma{G_2(F^{2}_{\tau} u) - G_1(F^{2}_{\tau}u)}_{\L1}
    \\
    & \leq &
    \O \cdot \left(
      \norma{Df_1 - Df_2}_{\C0(\Omega,\reali^{n\times n})}
      + 
      \norma{G_1-G_2}_{\C0(\mathcal{U}_{\delta_o};\L1(\reali;\reali^n))}
    \right) \,.
  \end{eqnarray*}
  Concerning the pointwise convergence, note that
  \begin{displaymath}
    \norma{F^{k}_t u - F_t u}_{\L1}
    \leq
    \mathcal{L} \int_0^t
    \norma{G(F_{\tau} u) - G_k(F_{\tau} u)}_{\L1}
    \;d\tau
  \end{displaymath}
  and the proof is concluded trough Lebesgue convergence theorem.
\end{proof}

\subsection{Integral Characterization}

Following~\cite[\S~9.2]{BressanLectureNotes},
\cite[\S~5.2]{AmadoriGosseGuerra}
and~\cite[Definition~15.1]{BianchiniBressan}, let $v \in
\mathcal{D}_T$. For all $\xi \in \reali$, define $U^\sharp_{(v,\xi)}$
as the solution to the homogeneous Riemann problem
\begin{equation}
  \label{eq:Sharp}
  \left\{
    \begin{array}{l}
      \partial_t w + \partial_x f(w) = 0
      \\
      w(0,x) = \left\{
        \begin{array}{l@{\quad\mbox{ if } \quad}rcl}
          \displaystyle \lim_{y \to \xi-} v(y) & x & < & \xi
          \\
          \displaystyle \lim_{y \to \xi+} v(y) & x & > & \xi \,.
        \end{array}
      \right.
    \end{array}
  \right.
\end{equation}
Define $U^\flat_{(v,\xi)}$ as the broad solution
(see~\cite[\S~3.1]{BressanLectureNotes}) to the Cauchy problem
\begin{equation}
  \label{eq:Flat}
  \left\{
    \begin{array}{l}
      \partial_t w + {Df \left( v(\xi) \right)} \, \partial_x w 
      = 
      G(v)
      \\
      w(0,x) = v(x) \,,
    \end{array}
  \right.
\end{equation}
so that denoting $l_i^{\xi} = l_i\left(v(\xi)\right)$,
$r_i^{\xi} = r_i\left(v(\xi)\right)$ and
$\lambda_i^{\xi} = \lambda_i \left( v(\xi)\right)$,
\begin{eqnarray*}
  U^\flat_{(v,\xi)} (t,x)
  & = &
  U^{\flat,1}_{(v,\xi)} (t,x) + U^{\flat,2}_{(v,\xi)} (t,x)
  \\
  U^{\flat,1}_{(v,\xi)} (t,x)
  & = &
  \sum_{i=1}^n  
  \left(
    l_i^{\xi} \cdot v\left(x-\lambda_i^{\xi} t\right)
  \right)
  r_i^{\xi}
  \\
  U^{\flat,2}_{(v,\xi)} (t,x)
  & = &
  \sum_{i=1}^n
  \int_{0}^{t} 
  \left(
    l_i^{\xi} \cdot 
    G (v) (x-\lambda_i^{\xi} s)
  \right)
  r_i^{\xi}\, ds \,.
\end{eqnarray*}
We are now ready to prove the first part of the characterization
stated in Theorem~\ref{thm:main}.

\begin{theorem}
  \label{thm:caratt1}
  Let $F$ be the map constructed in Theorem~\ref{thm:Conv}, let
  $\hat\lambda$ be an upper bound for all characteristic speeds. Then,
  for all $\bar t \in \left[0, T \right[$, $u \in \mathcal{D}_{\bar
    t}$ and all $\tau \in \left[0,T-\bar t\right]$, $F_\tau u$
  satisfies~\ref{it:thm:main} in Theorem~\ref{thm:main}.
\end{theorem}

\begin{proof}
  To obtain~\Ref{it:thm:main1}, compute:
  \begin{eqnarray*}
    & &
    \frac{1}{\theta} \,
    \int_{\xi - \theta \hat\lambda}^{\xi+ \theta \hat\lambda}
    \norma{ \left(F_\theta u(\tau)\right)(x) - U^\sharp_{(v,\xi)}(\theta,x)}
    \, dx
    \; \leq
    \\
    & \leq &
    \frac{1}{\theta} \,
    \int_{\xi - \theta \hat\lambda}^{\xi+ \theta \hat\lambda}
    \norma{ 
      \left(F_\theta u(\tau)\right) (x) - 
      \left(S_\theta u(\tau)\right) (x) -
      \theta\, \left(G\left( u(\tau) \right) \right) (x)
    }\, dx
    \\
    & &
    +
    \frac{1}{\theta} \,
    \int_{\xi - \theta \hat\lambda}^{\xi+ \theta \hat\lambda}
    \norma{
      \left(S_\theta u(\tau)\right) (x) -
      U^\sharp_{(v,\xi)}(\theta,x)} \, dx
    \\
    & &
    +
    \frac{1}{\theta} \,
    \int_{\xi - \theta \hat\lambda}^{\xi+ \theta \hat\lambda}
    \norma{\theta\, \left(G\left( u(\tau) \right) \right) (x)} \, dx
    \\
    & \leq &
    \frac{1}{\theta} \,
    \norma{ 
      F_\theta u(\tau) - 
      S_\theta u(\tau) -
      \theta\, G\left( u(\tau) \right)
    }_{\L1}
    \\
    & &
    +
    \frac{1}{\theta} \,
    \int_{\xi - \theta \hat\lambda}^{\xi+ \theta \hat\lambda}
    \norma{
      \left(S_\theta u(\tau)\right) (x) -
      U^\sharp_{(v,\xi)}(\theta,x)} \, dx
    \\
    & &
    +
    \int_{\xi - \theta \hat\lambda}^{\xi+ \theta \hat\lambda}
    \norma{\left(G\left( u(\tau) \right) \right) (x)} \, dx
  \end{eqnarray*}
  As $\theta \to 0$, the first summand above vanishes by~\Ref{it:tg}
  in Theorem~\ref{thm:main}, the second by~\cite[(9.16),
  \S~9.2]{BressanLectureNotes} and the latter one by
  $G\left(u(\tau)\right) \in \L1$.

  Similarly, to obtain~\Ref{it:thm:main2}, we
  exploit~\cite[\S~9.2]{BressanLectureNotes}:
  \begin{eqnarray*}
    & &
    \frac{1}{\theta} \,
    \int_{a + \theta\hat\lambda}^{b - \theta \hat\lambda} 
    \norma{\left(F_\theta u(\tau)\right)(x) - 
      U^\flat_{(v,\xi)} (\theta,x)} \, dx    
    \quad \leq
    \\
    & \leq &
    \frac{1}{\theta} \,
    \int_{a + \theta\hat\lambda}^{b - \theta \hat\lambda}
    \norma{ 
      \left(F_\theta u(\tau)\right) (x) - 
      \left(S_\theta u(\tau)\right) (x) -
      \theta\, \left(G\left( u(\tau) \right) \right) (x)
    }\, dx
    \\
    & &
    +
    \frac{1}{\theta} \,
    \int_{a + \theta\hat\lambda}^{b - \theta \hat\lambda}
    \norma{
      \left(S_\theta u(\tau)\right) (x) -
      U^\flat_{(v,\xi)}(\theta,x)
      +
      \theta\, \left(G\left( u(\tau) \right) \right) (x)} \, dx
    \\
    & \leq &
    \frac{1}{\theta} \,
    \norma{ 
      F_\theta u(\tau) - 
      S_\theta u(\tau) -
      \theta\, G\left( u(\tau) \right)
    }_{\L1}
    \\
    & &
    +
    \frac{1}{\theta} \,
    \int_{a + \theta\hat\lambda}^{b - \theta \hat\lambda}
    \norma{
      \left(S_\theta u(\tau)\right) (x) -
      U^{\flat,1}_{(v,\xi)}(\theta,x)} \, dx
    \\
    & &
    +
    \frac{1}{\theta} \,
    \int_{a + \theta\hat\lambda}^{b - \theta \hat\lambda}
    \norma{
      \theta \left(G\left( u(\tau) \right) \right) (x) -
      U^{\flat,2}_{(v,\xi)}(\theta,x)} \, dx
    \\
    & \leq &
    \frac{1}{\theta} \,
    \norma{ 
      F_\theta u(\tau) - S_\theta u(\tau) - \theta\, G\left( u(\tau) \right)
    }_{\L1}
    \\
    & &
    +
    \O \left( \tv \left( u(\tau) ; \left]a,b\right[\right) \right)^2
    \\
    & &
    +
    \frac{1}{\theta} \,
    \int_{a + \theta\hat\lambda}^{b - \theta \hat\lambda}
    \norma{
      \theta \left(G\left( u(\tau) \right) \right) (x) -
      U^{\flat,2}_{(v,\xi)}(\theta,x)} \, dx
  \end{eqnarray*}
  As $\theta \to 0$, the first summand above vanishes by~~\Ref{it:tg}
  in Theorem~\ref{thm:main}, while the latter vanishes as $\theta \to
  0$. Indeed, using~\cite[Lemma~2.3]{BressanLectureNotes},
  \begin{eqnarray*}
    & & 
    \frac{1}{\theta} \,
    \int_{a + \theta\hat\lambda}^{b - \theta \hat\lambda}
    \norma{
      \theta \left(G\left( u(\tau) \right) \right) (x) -
      U^{\flat,2}_{(v,\xi)}(\theta,x)} \, dx
    \leq
    \\
    & \leq &
    \frac{1}{\theta} \,
    \int_{a + \theta\hat\lambda}^{b - \theta \hat\lambda}
    \Bigg\Vert 
    \int_0^\theta 
    \sum_{i=1}^n 
    l_i^{(\tau,\xi)} \cdot \left( G \left( u(\tau) \right) (x) \right)
    r_i^{(\tau,\xi)} \, ds
    \\
    &&
    \qquad\qquad
    -
    \sum_{i=1}^n
    \int_0^\theta 
    l_i^{(\tau,\xi)} \cdot 
    \left( G \left( u(\tau) \right) \right)(x- \lambda_i^{(\tau,\xi)}s)
    r_i^{(\tau,\xi)} \, ds\Bigg\Vert dx
    \\
    & \leq &
    \frac{\O}{\theta} \sum_{i=1}^n \int_0^\theta \!
    \int_{a+\theta\hat\lambda}^{b-\theta\hat\lambda}
    \norma{
      \left( G \left( u(\tau) \right) \right) (x) - 
      \left( G \left( u(\tau) \right) \right) (x-\lambda_i^{(\tau,\xi)}s)
    }
    \, dx \, ds
    \\
    & \leq &
    \frac{\O}{\theta} \sum_{i=1}^n \int_0^\theta
    \tv \left( G \left( u(\tau) \right); \left] a,b\right[ \right)
    \modulo{\lambda_i^{(\tau,\xi)}} \, s \, ds
    \\
    & \leq &
    \O \cdot \theta \cdot
    \tv \left( G \left( u(\tau) \right); \left] a,b\right[ \right) \,.
  \end{eqnarray*}
  completing the proof of~\Ref{it:thm:main2}.
\end{proof}

\begin{theorem}
  \label{thm:caratt2}
  Let $F$ be the map constructed in Theorem~\ref{thm:Conv} and
  $\hat\lambda$ be an upper bound for all characteristic speeds.
  Assume $\bar t \in \left[0, T \right[$, $u \colon [0,T-\bar t]
  \mapsto \mathcal{D}_{T}$ is Lipschitz, satisfies $u(t) \in
  \mathcal{D}_{\bar t +t}$ and both~\Ref{it:thm2:main1}
  and~\Ref{it:thm2:main2} in Theorem~\ref{thm:main} hold. Then, $u(t)
  = F_t u$ for all $t \in [0,T-\bar t]$.
\end{theorem}

We omit this proof, since it is a slight modification of~\cite[Part~2
of Theorem~9.2]{BressanLectureNotes}.

\subsection{Consequences of the Hyperbolic Rescaling}

Given a function $v \colon \reali \mapsto \reali^n$ and $\lambda >0$,
we denote by $v_\lambda$ the function obtained by applying a
dilatation to $v$, i.e. $v_\lambda(x)=v(\lambda x)$. Obviously
$v\in\overline{\mathcal{D}}_\delta$ implies
$v_\lambda\in\overline{\mathcal{D}}_\delta$.  We have the following
Proposition (see also~\cite[Corollary~1]{MR776388}).

\begin{proposition}
  \label{prop:Guerra}
  Let $S\colon [0,T] \times \overline{\mathcal{D}}_\delta \mapsto
  \overline{\mathcal{D}}_\delta$ be the semigroup generated by a
  system of conservation laws and let $d\colon
  \overline{\mathcal{D}}_\delta \times \overline{\mathcal{D}}_\delta
  \mapsto \reali^+$ be a distance satisfying
  \begin{displaymath}
    d(u_\lambda, v_\lambda)
    =
    \frac{1}{\lambda} \, d(u,v)
    \mbox{ for all } u,v\in\overline{\mathcal{D}}_\delta,
    \mbox{ and }\lambda>0,
  \end{displaymath}
  and the Gr\"onwall estimate for a positive $C$:
  \begin{displaymath}
    d(S_tu,S_tv) \leq e^{Ct} \, d(u,v)
    \mbox{ for all } u,v \in \overline{\mathcal{D}}_\delta
    \mbox{ and } t\in [0,T].
  \end{displaymath}
  Then $C=0$, i.e.~$d$ is non expansive with respect to $S$.
\end{proposition}

\begin{proof}
  If $u(t,x) = \left(S_t u\right)(x)$ is a semigroup trajectory, then
  also $u(\lambda t,\lambda x) = \left(S_{t} u_\lambda\right)(x)$ is a
  semigroup trajectory. Therefore we have the equality
  \begin{displaymath}
    \left(S_tu\right)_\lambda  = S_{\frac{t}{\lambda}} u_\lambda
    \mbox{ for all }
    u\in\overline{\mathcal{D}}_\delta,\ t\in[0,T]
    \mbox{ and }
    \lambda>0 \,.
  \end{displaymath}
  Hence we can compute for all $u,v \in \overline{\mathcal{D}}_\delta$
  and $\lambda > 0$
  \begin{displaymath}
    \begin{array}{rcccl}
      d\left(S_tu,S_tv\right)
      & = &
      \lambda \, 
      d\left( \left(S_tu\right)_\lambda, \left(S_tv\right)_\lambda\right)
      & = &
      \lambda \,
      d\left( S_{\frac{t}{\lambda}} u_\lambda, S_{\frac{t}{\lambda}}v_\lambda \right)
      \\
      &\leq &
      \lambda \, e^{C\frac{t}{\lambda}} d\left( u_\lambda,v_\lambda \right)
      & = &
      e^{C\frac{t}{\lambda}}d\left(u,v\right)
    \end{array}
  \end{displaymath}
  Now, letting $\lambda$ tend to infinity, we get the non expansive
  property
  \begin{displaymath}
    d\left(S_t u,S_t v\right)
    \leq 
    d\left(u,v\right)
    \mbox{ for all } u,v \in \overline{\mathcal{D}}_\delta
    \mbox{ and } t \in [0,T].
  \end{displaymath}
\end{proof}

{\small{

    \bibliographystyle{abbrv}

  }}

\end{document}